\documentclass[12pt]{article}
\usepackage{amsmath,amssymb,amsfonts,chngcntr}

\usepackage{rsfso}
\usepackage{framed}
\newcommand{\dlines}{\displaylines}
\def\limn{\lim_{n\to\infty}}

\def\paref#1{(\ref{#1})}

\def\eqref#1{(\ref{#1})}

\newcommand{\R}{\mathbb{R}}
\renewcommand{\th}{\theta}

\newcommand{\N}{\mathbb{N}}

\newcommand{\Var}{{\rm Var}}
\newcommand{\finedim}{\par\hfill$\blacksquare$\hphantom{aaaaaa}\par\noindent\ignorespaces}
\let\qed=\finedim

\let\dstyle=\displaystyle

\newcommand{\e}{{\rm e}}

\def\tin#1{\par\noindent\hskip3em\llap{#1\enspace}\ignorespaces}
\def\cl#1{{\mathcal #1}}

\def\E{{\rm E}}

\def\P{{\rm P}}

\def\tfrac#1#2{{\textstyle\frac {#1}{#2}}}
\newtheorem{theorem}{Theorem}[section]
\newtheorem{rema}{Remark}[section]

\newtheorem{cor}{Corollary}[section]

\newtheorem{example}{Example}[section]
\newtheorem{lem}{Lemma}[section]

\newcommand{\proof}{{\it Proof}.\ \/}
\makeatletter
\@addtoreset{equation}{section}
\makeatother

\begin{document}
\title{\bf Tightness and exponential tightness of Gaussian probabilities}
\author{Paolo Baldi\protect\footnote{baldi@mat.uniroma2.it, http://mat.uniroma2.it/$\sim$processi/baldi.htm}
\protect\footnote{The author acknowledges the MIUR Excellence Department Project awarded to the Dipartimento di Matematica, Universit\`a  di  Roma  ``Tor  Vergata'',  CUP  E83C18000100006
}\\
{\it Dipartimento di Matematica, Universit\`a
di Roma Tor Vergata, Italy}\\
} \date{}\maketitle
\begin{abstract}
We prove a simple criterion of exponential tightness for sequences of Gaussian r.v.'s with values in a separable Banach space from which
we deduce a general result of Large Deviations which allows easily to obtain LD estimates in various situations.
\end{abstract}

\smallskip


\section{Introduction}
Let $(X_n)_n $ be a sequence of $E$-valued Gaussian r.v.'s, where $E$
is a separable Banach space. This article is concerned with the investigation
of Large Deviation properties at speed $g$ of the sequence
$(g_n^{-1/2}X_n)_n $, where $g:\N\to\R^+$ is a speed function i.e. such that $\lim_{n\to
\infty}g_n=+\infty$.

This is a classical problem that is completely understood if the law
of $X_n $ does not depend on $n$ (see \cite{dvIII} \S 5 or the
exposition in \cite{az-stflour}, p.39).

In \cite{ampa} (see also \cite{DZ10} \S4.5.3) it is proved that if, for every $\th\in E'$,
\begin{equation}\label{log-Laplace}
\lim_{n\to \infty} \frac 1{g_n}\log\E[\e^{g_n\langle\th,Y_n
\rangle}]:=\Lambda(\th)
\end{equation}
{\it and} the sequence $(Y_n)_n $ is exponentially tight,
then it enjoys some Large Deviations estimates at speed $g$ with a
rate function $I$ which is the convex conjugate of $\Lambda$
(see the precise statement in Theorem \ref{egb} below).

It is therefore of interest to produce criteria ensuring the
property of exponential tightness.

Recall that a sequence of r.v.'s $(X_n)_n $ with values in some metric
space $E$ is said to be {\it tight} as $n\to +\infty$ if for every $\delta>0$
there exists a compact set $K_\delta$ such that
\begin{equation}\label{et0}
\P(X_n \not\in K_\delta)<\delta\ .
\end{equation}
A sequence of r.v.'s $(X_n)_n $ with values in some metric space is
said to be {\it exponentially tight} as $ n \to \infty$ at speed $g$ if
for every $R>0$ there exists a compact set $K_R$ such that
\begin{equation}\label{et-exp}
\P(X_n \not\in K_R)<\e^{-g_nR}\ .
\end{equation}
The main result of this paper (Theorem \ref{mt}) is that if a sequence
$(X_n)_n$ of Gaussian $E$-valued r.v.'s is tight, then, for every
speed function $g$, the sequence $(g_n^{-1/2}X_n)_n $ is
exponentially tight at speed $g$.

This result will allow us to prove the following result.
\begin{theorem}\label{mt2}
Let $(X_n)_n $ be a sequence of Gaussian r.v.'s with values in the
separable Banach space $E$ and converging in law to a r.v. $X$.
Then, for every speed function $g$, the sequence
$(g_n^{-1/2}X_n)_n $ satisfies a Large Deviation Principle with
speed $g$ and rate function
\begin{equation}\label{cc}
I(x)=\frac 12\Vert x\Vert^2_{\cl H}
\end{equation}
where $\Vert \enspace\Vert_{\cl H}$ denotes the Replicating Kernel Hilbert Space norm of the law of the r.v. $X-\E[X]$, with
the understanding
$\Vert x\Vert_{\cl H}=+\infty$ if $x\not\in\cl H$.
\end{theorem}
In some sense Theorem \ref{mt2} states that if $X_n\to X$ in law, then, in a Gaussian setting, the Large Deviations
asymptotics of
the sequences $(g_n^{-1/2}X_n)_n$ and $(g_n^{-1/2}X)_n $ are the same ({\it without any constraint}
concerning the speed of convergence of $X_n\to X$).

Theorem \ref{mt2} is well known if the law
of $X_n $ does not depend on $ n$ (see \cite{dvIII} \S 5 or the
exposition in \cite{az-stflour}, p.39). This, of course, suggests
that Theorem \ref{mt2} should be an immediate consequence of the
above mentioned result for the sequence $(g_n^{-1/2}X)_n $, using
the argument of exponential approximation, as explained in \cite{DZ10}
\S4.2.2. Actually it turns out
that Theorem \ref{mt2} can be proved using this argument, see
\S\ref{exp-appr}, but this way of proof would also  require the exponential tightness
result of Theorem \ref{mt}.

Recall that Fernique's theorem states that a Gaussian r.v. has
bounded exponential moments. As a, possibly useful, by-product of
our treatment, we obtain (Corollary \ref{f-cor}) that a tight family
of $E$-valued Gaussian r.v.'s has uniformly bounded exponential
moments.

In \S\ref{et-main} we prove our main result (Theorem \ref{mt}) i.e. that if $(X_n)_n $ is tight then
$(g_n^{-1/2}X_n)_n $ is exponentially tight at speed $g$.

In \S\ref{appl} and \S\ref{exp-appr} we give two different proofs of Theorem \ref{mt2} whereas section
\ref{concl} is devoted to examples and counterexamples.
\section{Exponential tightness of Gaussian families}\label{et-main}
Throughout this paper $E$ shall denote a separable Banach space,
$E'$ its topological dual.

In order to investigate exponential tightness of Gaussian families we shall take advantage of
the following result which is a particular case of Theorem 3.1 of De Acosta \cite{deaco}.
\begin{theorem}\label{deacoth}{\rm (\cite{deaco})} Let $(\mu_\alpha)_{\alpha\in {\cl A}}$ be a family
of probabilities on the separable Banach space $E$ such that
\tin{a)}$(\mu_\alpha)_{\alpha\in {\cl A}}$ is tight.
\tin{b)} There exists $t>0$ such that
\begin{equation}\label{cond-deaco}
\sup_{\alpha\in {\cl A}}\int_E \e^{t|x|}\, d\mu_\alpha(x)<+\infty\ .
\end{equation}
Then there exists a convex compact well balanced set $K\subset E$ such that
$$
\sup_{\alpha\in {\cl A}}\int_E \e^{q_K(x)}\, d\mu_\alpha(x)<+\infty
$$
where $q_K$ denotes the Minkowski functional of the set $K$.
\end{theorem}
The definition of a Minkowski functional will be recalled shortly.

Remark that in Theorem \ref{deacoth} there is no assumption of Gaussianity. In order to take advantage of Theorem \ref{deacoth},
we first prove that condition \eqref{cond-deaco} is automatically
satisfied if Gaussianity, in addition to the tightness condition of assumption a) of the theorem, is enforced.

This fact will follow from the next result, whose proof follows the
same line of reasoning of the classical Fernique's theorem (see
\cite{fernique-cras}, \cite{fernique-stflour74} p. 11).
\begin{theorem} \label{fernique}
Let $E$ be a separable Banach space and $\phi: E\to[0,+\infty]$ a
measurable semi-norm, i.e. a measurable application such that
\begin{equation}\label{semin}
\begin{array}{rl}
 \phi(\lambda x)&= |\lambda|\phi(x),\quad \lambda\in\R\cr
\phi(x+y)&\le\phi(x)+\phi(y)\ .
\end{array}
\end{equation}
Let $(\mu_\alpha)_{\alpha\in {\cl A}}$ be a family of centered Gaussian probabilities
on $E$ such that $\mu_\alpha(\phi<+\infty)=1$ for every
$\alpha\in \cl A$. If there exists $s>0$ such that
$\mu_\alpha(\phi>s)\le\beta<\frac12$ for every $\alpha\in \cl A$, then
there exists $a>0$ such that
$$
\sup_{\alpha\in \cl A}\int_E \e^{ a \phi^2(x)}\,\mu_\alpha(dx)<\infty\ .
$$
\end{theorem}
\proof From the relation
\begin{equation}\label{ps-tightness}
\mu_\alpha(\phi> s)\le\beta<\frac 12\ ,
\end{equation}
following the lines of the proof of Fernique's theorem (see also
\cite{deuschel-stroock} \S1.3) if we define by recurrence the sequence $(t_n)_n$ by
$$
t_0=s\qquad t_n=\sqrt{2}\,t_{n-1}+s
$$
i.e.
\begin{equation}\label{bzeta}
t_n=s\bigl(1+\sqrt{2}+\cdots+\bigl(\sqrt{2}\bigr)^n\bigr)=
s\frac{(\sqrt{2})^{n+1}-1}{\sqrt{2}- 1}\leq
\underbrace{s\Bigl(\frac{\sqrt{2}}{\sqrt{2}-1}\Bigr)}_{\displaystyle
:=\sqrt{\zeta}} 2^{n/2}\ .
\end{equation}
we obtain, for every $\alpha\in \cl A$,
\begin{equation}\label{sch.227}
\mu_\alpha(\phi> t_n)\le\kappa^{2^n}\mu_\alpha(\phi< s)
\end{equation}
where
$$
\kappa=\frac\beta{1-\beta}<1
$$
We can now split the integral
$$
\int_E \e^{ a
\phi^2(x)}\,\mu_\alpha(dx)\leq\underbrace{\int_{\{\phi^2\leq
\zeta\}} \e^{ a \phi^2(x)}\,\mu_\alpha(dx)}_{\le\e^{ a
\zeta}}+\int_{\{\phi^2> \zeta\}}\e^{ a \phi^2(x)}\,\mu_\alpha(dx)\ .
$$
But
$$
\dlines{\int_{\{\phi^2> \zeta\}}\e^{ a
\phi^2(x)}\,\mu_\alpha(dx)\leq\sum_{n=0}^\infty\int_{\{\zeta\cdot
2^n<\phi^2\leq \zeta\cdot 2^{n+1}\}}\e^{ a
\phi^2(x)}\,\mu_\alpha(dx)\leq\cr \le\sum_{n=0}^\infty\e^{ a
\zeta\cdot 2^{n+1}}\mu_\alpha(\phi^2> \zeta\cdot 2^n)\ .\cr }
$$
Since $t_n^2\leq \zeta 2^n$, \eqref{sch.227} gives
$$
\dlines{\int_E\e^{ a \phi^2(x)}\,\mu_\alpha(dx)\leq \e^{a\zeta}+
\sum_{n=0}^\infty\e^{ a\zeta\cdot 2^{n+1}}\mu_\alpha(\phi> t_n)
\le \e^{ a \zeta}+\sum_{n=0}^\infty\e^{ a  \zeta\cdot
2^{n+1}}\mu_\alpha(\phi< s)\cdot \kappa^{2^n}\cr }
$$
which, for $ a  < \frac {\log \frac1\kappa}{2\zeta}$, gives a convergent  series thus
concluding the proof.
\finedim
\begin{cor}\label{f-cor} Let $(\mu_\alpha)_{\alpha\in {\cl A}}$ be a family of centered Gaussian probabilities
on the Banach space $E$, such that there exists $s>0$ such that $\mu_\alpha(|x|>s)\le \beta<\frac12$
for every $\alpha\in {\cl A}$. Then
\begin{equation}\label{maggio}
\sup_{\alpha\in {\cl A}} \int_E \e^{t|x|}\, d\mu_\alpha(x)<+\infty,\qquad \mbox{for every }t>0\ .
\end{equation}
In particular \eqref{maggio} holds if the family $(\mu_\alpha)_{\alpha\in {\cl A}}$ is tight.
\end{cor}
\proof Let $a>$0 be such that
$$
\sup_{\alpha\in {\cl A}} \int_E \e^{a|x|^2}\, d\mu_\alpha(x)<+\infty
$$
as guaranteed by Theorem \ref{fernique}, then
$$
\dlines{
\int_E \e^{t|x|}\, d\mu_\alpha(x)= \int_{\{|x|\le t/a\}} \e^{t|x|}\, d\mu_\alpha(x)+
\int_{\{|x|> t/a\}} \e^{t|x|}\, d\mu_\alpha(x)\le\cr
\le \e^{t^2/a}+\int_E \e^{a|x|^2}\, d\mu_\alpha(x)<+\infty\ .\cr
}
$$
Finally remark that the condition $\mu_\alpha(|x|>s)\le \beta<\frac12$ for every $\alpha\in{\cl A}$ is automatically satisfied if
$(\mu_\alpha)_{\alpha\in\cl A}$ is tight, as compact sets are bounded in $E$.
\finedim
Let $K\subset E$. Let us recall that its Minkowski functional is a map
$q_K:E\to [0,+\infty]$ defined as
$$
q_K(x)=\inf\{t>0; x\in tK\}\ .
$$
It is well known that if $K$ is a convex set then $q_K$ is subadditive and that if $K$ is well balanced then
$q_K$ is positively homogeneous. Well balanced means that if $x\in K$ then also $t x\in K$ for every $0\le t\le 1$.

The following is the main result of this section.
\begin{theorem}\label{mt}
Let $(X_\alpha)_{\alpha\in \cl A}$ be a tight family of Gaussian $E$-valued r.v.'s. Then for every speed function $g$ the family
$(g_n^{-1/2}X_\alpha)_n$ is uniformly exponentially tight, i.e., for every $R>0$ there exists a compact set $K_R\subset E$ such that
$$
\P(g_n^{-1/2}X_\alpha\not\in K_R)\le \e^{-g_n R}
$$
for every $\alpha\in \cl A$ and for every $n$. In particular, if a sequence $(X_n)_n$ of Gaussian $E$-valued r.v.'s is tight, then
for every speed function $g$, $(g_n^{-1/2}X_n)_n$  is exponentially tight.
\end{theorem}
\proof
Let us assume first that the r.v.'s $X_\alpha$ are centered.
By Corollary \ref{f-cor} also condition b) of Theorem
\ref{deacoth} is satisfied. Thus, by Theorem \ref{deacoth} applied to the family $(\mu_\alpha)_{\alpha\in\cl  A}$
of the laws of the $X_\alpha$'s, there exists a  convex compact well
balanced set $K$ such that
$$
\sup_{\alpha\in\cl  A}\E\bigl[\e^{q_K(X_\alpha)}\bigr]:=C<+\infty\ .
$$
For every $\alpha\in\cl  A$ let $(Z_n)_n$ be a sequence of i.i.d. r.v.'s with the same law as $X_\alpha$. Let us assume
for simplicity, at first,
that $g_n$ is an integer number for every $n$.
We have
\begin{equation}\label{intg}
g_n^{-1/2}(Z_1+\dots +Z_{g_n})\sim X_\alpha
\end{equation}
hence 
$$
\dlines{
\E\bigl[\e^{q_K(g_n^{1/2}X_\alpha)}\bigr]=\E\bigl[\e^{q_K(Z_1+\dots +Z_{g_n})}\bigr]\le
\E\bigl[\e^{q_K(Z_1)+\dots+q_K(Z_{g_n})}\bigr]=
=\bigl(\E\bigl[\e^{q_K(X_\alpha)}\bigr]\bigr)^{g_n}\le C^{g_n}\ .\cr
}
$$
Remark the relation
\begin{equation}\label{mkw1}
\{x\not\in t K\}=\{x;q_K(x)>t\}\ .
\end{equation}
Hence for every $t>0$, by Markov inequality,
$$
\dlines{ \P\bigl(g_n^{-1/2} X_\alpha\not\in t K\bigr)=\P\bigl(g_n^{1/2}
X_\alpha\not\in t g_n K\bigr)= \P\bigl(q_K(g_n^{1/2} X_\alpha)\ge t
g_n\bigr)\le\cr \le\e^{-t g_n}\E\bigl[\e^{q_K(g_n^{1/2}X_\alpha)}\bigr]
\le \e^{(\log C-t) g_n}\cr }
$$
and if we chose $t=R+\log C$, the compact $K_R=t K$ satisfies the
requirement \eqref{et-exp}.

In general, if $g$ is not integer valued, \eqref{intg} will become
$g_n^{-1/2}(Z_1+\dots+Z_{\lfloor g_n\rfloor})\sim \lfloor g_n\rfloor^{1/2} g_n^{-1/2}X_\alpha$, $\lfloor\enspace\rfloor$
denoting the
integer part function and the remainder of the proof is to be modified accordingly.

If we do not assume the $X_\alpha$'s to be centered, let us denote by $x_\alpha=\E[X_\alpha]$ their means and let
$\widetilde X_\alpha=X_\alpha-x_\alpha$. The assumption of tightness implies that
the sequence $(x_\alpha)_\alpha$ is contained in some compact set $K$, that we can assume to be well balanced,
and that also the sequence
$(\widetilde X_\alpha)_\alpha$ is tight. Let $\widetilde K_R$ be a compact, convex and well balanced set
such that $\P(g_n^{-1/2}\widetilde X_\alpha\not\in \widetilde K_R)\le \e^{-g_n R}$.
Then $K_R=\widetilde K_R+K$ is a compact set and we have
$$
\P\bigl(g_n^{-1/2} X_\alpha\not\in K_R\bigr)\le \P\bigl(g_n^{-1/2} \widetilde X_\alpha\not\in \widetilde K_R\bigr)\le\e^{-g_nR}\ .
$$
\finedim
It is natural, at this time, to inquire whether a converse of Theorem \ref{mt} holds. It is obvious that uniform exponential tightness, even with respect to a single speed function implies tightness.

Whether a Gaussian sequence $(X_n)_n$ that is exponentially tight with respect to every speed function is also tight is less obvious, as exponential tightness of the sequence is a weaker property than the uniform exponential tightness of the family $(X_n)_n$. A partial answer to this question is given in Example \ref{counter}.
\section{Large Deviations}\label{appl}%
We give two proofs of Theorem \ref{mt2}. The second one, in the next section, is much shorter but maybe the first one gives more insight
into the structure of Gaussian probabilities. Both proofs rely on Theorems \ref{egb} and \ref{mt}.
\medskip

In order to prove Theorem \ref{mt2} we shall take advantage of the following result.

Let us recall that $x\in E$  is said to be an {\it exposed point} for the rate function
$I$ if there exists $\eta\in E'$ such that, for $x\not =y$,
\begin{equation}\label{eg.97}
I(y)-I(x)>\langle\eta,y- x\rangle
\end{equation}
which is a condition of strict convexity of the rate function $I$.
$\eta$ is called an {\it exposing hyperplane}. The following theorem is an infinite dimensional extension of the
Ellis-G\"artner theorem (\cite{gart}, \cite{ellis}).
\begin{theorem} \label{egb} {\rm (\cite{ampa}, see also \cite{DZ10} \S4.5.3)} Let $(\mu_n)_n $ be an exponentially
tight sequence of probabilities at speed $g$ on the separable Banach space
$E$, such that, for every $\th\in E'$, the limit
\begin{equation}\label{eg.177}
\lim_{ n \to\infty}\frac 1{g_n}\log\int_E\e^{g_n\langle\th,x\rangle}\,d\mu_n(x):=
\Lambda(\th)
\end{equation}
exists. Let us denote
\begin{equation}\label{rate1}
I(x)=\sup_{\th\in E'}(\langle \th,x\rangle-\Lambda(\th))
\end{equation}
the convex conjugate of $\Lambda$.
Then for every closed set $F\subset E$
\begin{equation}\label{eg.up0}
\liminf_{ n \to\infty} \frac 1{g_n}\log\mu_n(F)\le -\inf_{x\in F}I(x)
\end{equation}
and for every open set $G\subset X$,
\begin{equation}\label{eg.low0}
\liminf_{ n \to\infty} \frac 1{g_n}\log\mu_n(G)\ge -\inf_{x\in G\cap\cl E}I(x)
\end{equation}
where $\cl E$ denotes the set of all exposed points of $I$
such that for an exposing hyperplane $\eta$ it holds
$\Lambda(s\eta)<+\infty$ for some $s>1$.
\end{theorem}
The application of Theorem \ref{egb}, in order to obtain a Large Deviation Principle, is not always undemanding because
it provides an incomplete lower bound. We shall see however that from \eqref{eg.low0}
a true lower LDP bound can be easily obtained if the probabilities $\mu_n$ are Gaussian.
\bigskip

\noindent{\it First proof} of Theorem \ref{mt2}.
By Theorem \ref{mt} the sequence $(g_n^{-1/2}X_n)_n$ is exponentially tight at speed $g$.
Let us denote $x_n=\E[X_n]$, $\widetilde X_n=X_n-x_n$, $x=\E[X]$, $\widetilde X=X-x$. We have
$$
\displaylines{
\log\E[\e^{g_n\langle \th,g_n^{-1/2}X_n \rangle}]=
g_n^{1/2}\langle \th, x_n\rangle+\log\E[\e^{g_n\langle \th,g_n^{-1/2}\widetilde X_n \rangle}]\cr
}
$$
Hence condition (\ref{eg.177}) is satisfied as the r.v.'s
$\langle \th,X_n \rangle$, $\th\in E'$, are Gaussian so that
\begin{equation}\label{Lambda}
\begin{array}{c}
\dstyle\lim_{ n \to\infty}\frac 1{g_n}\,\log\E[\e^{g_n\langle\th,g_n^{-1/2}X_n \rangle}]=\lim_{ n
\to\infty}\frac 1{g_n}\Bigl(g_n^{1/2}\langle \th, x_n\rangle+\log\E[\e^{g_n\langle \th,g_n^{-1/2}\widetilde X_n \rangle}]\Bigr)=\cr
=\dstyle\lim_{ n \to\infty}\frac 1{2g_n}\Var(g_n^{1/2}\langle \th,\widetilde X_n\rangle)=\frac 12\,\Var(\langle \th,\widetilde X\rangle)=\log\E[\e^{\langle\th, \widetilde X\rangle}]:=\Lambda(\th)\ .
\end{array}
\end{equation}
In order to conclude the proof of Theorem \ref{mt2} we must compute the convex conjugate $I$ of $\Lambda$.
This is well known (see \cite{az-stflour} Chap II e.g.), however, in order to investigate the
exposed points of $I$ we recall some facts about this computation and about
Gaussian probabilities $\mu$ on a separable Banach space $E$.

For every $\th\in E'$, the function $E\ni z\mapsto \langle \th, z\rangle$
defines a r.v. on the probability space $(E, \mu)$.
Hence $E'$ can be considered as a Gaussian space of r.v.'s on $(E, \mu)$.
Let us denote by $\overline{E}'_\mu$ the closure of $E'$ in the Hilbert space $L^2(\mu)$.
$\overline{E}'_\mu$ is a space of Gaussian r.v.'s and for every $\phi\in  \overline{E}'_\mu$ we can define its barycenter
\begin{equation}\label{by}
h=\int_E z\phi(z)\, d\mu(z)\ .
\end{equation}
The convergence of the integral in (\ref{by}) is a consequence of well known facts about the existence of moments of Gaussian probabilities.
The elements $h\in E$ of this form constitute a vector space $\cl H$ that becomes an Hilbert space when endowed with the norm
\begin{equation}\label{gd}
\Vert h\Vert_{\cl H}=\Vert\phi\Vert_{L^2(\mu)}=\int_E\phi(z)^2\,d\mu(z)
\end{equation}
which makes $\cl H$ isometric with $\overline{E}'_\mu$. Actually it is easy to see that \eqref{gd} is a good definition:
if $h$ is the barycenter
of $\phi$ and of $\phi'$, $\phi,\phi'\in \overline{E}'_\mu$, then $\Vert\phi\Vert_{L^2(\mu)}=\Vert\phi'\Vert_{L^2(\mu)}$ necessarily.

It is important to remark that the Hilbert space $\cl H$ is dense in the closed
subspace $\widetilde E=\mathop{\rm supp}\mu$. Actually if $\th\in E'$ is such that $\langle\th,h\rangle=0$ for every $h\in\cl H$,
then, if $h$
is the barycenter of $\phi\in \overline{E}'_\mu$, we would have
$$
0=\langle\th,h\rangle=\int_E\langle\th,x\rangle\phi(x)\,d\mu(x)\ .
$$
The function $x\mapsto\langle\th,x\rangle$ being orthogonal to every $\phi\in \overline{E}'_\mu$, necessarily $\langle\th,x\rangle=0$
$\mu$-a.s.

The functional $\Lambda:E'\to\R $ defined in \eqref{Lambda} is convex and,
as we already remarked, the computation of its convex conjugate $I$, i.e.
\begin{equation}\label{attained}
I(x)=\sup_{\th\in E'}(\langle\th,x\rangle-\Lambda(\th))
\end{equation}
is a classical fact (see \cite{az-stflour}, Proposition 1.5 p.53 and the literature therein).
Actually, if $x\in\cl H$ and $x=\int z\phi(z)\, d\mu(z)$ for some $\phi\in\overline{E}'_\mu$, we have
\begin{equation}\label{attained2}
\begin{array}{c}
\dstyle I(x)=\sup_{\th\in E'}\Bigl(\int_E\langle\th,z\rangle\phi(z)\,d\mu(z)-\frac 12\int_E\langle\th,z\rangle^2\,d\mu(z)\Bigr)=\\
\dstyle=\sup_{\th\in E'}\Bigl(-\frac12\int(\langle\th,z\rangle-\phi(z))^2\,d\mu(z)\Bigr)+\frac 12\int_E\phi(z)^2\,d\mu(z)
\end{array}
\end{equation}
from which, by considering a sequence $(\th_n)_n\subset E'$ converging to $\phi$ in $\overline{E'_\mu}$, we have
$$
\sup_{\th\in E'}\Bigl(-\frac12\int(\langle\th,z\rangle-\phi(z))^2\,d\mu(z)\Bigr)=0
$$
and
$$
I(x)=\frac 12\int_E\phi(z)^2\,d\mu(z)=\frac 12\Vert x\Vert_{\cl H}^2\ .
$$
Moreover if $\phi\in E'$ (which is a stronger assumption than $\phi\in \overline{E}'_\mu$) the supremum in \eqref{attained2} is attained at $\th=\phi$ so that
in this case we have also $I(x)=\langle\phi,x\rangle-\Lambda(\phi)=\Lambda(\phi)$. See \cite{az-stflour} Proposition 1.5 p.~53 for a proof that
$I(x)=+\infty$ if $x\not\in\cl H$. Then, thanks to Theorem \ref{egb},
the Large Deviation upper bound \eqref{eg.up0} holds with respect to the rate function \eqref{cc}.
In order to complete the proof of Theorem \ref{mt2} we must show that the exposed points, $\cl E$, of $I$ satisfying
the condition of Theorem \ref{egb} are such that, for every open set $A\subset E$,
$$
\inf_{x\in {A}\cap\cl E}I(x)=\inf_{x\in {A}}I(x)\ .
$$
This is the consequence of the following lemma.
\qed
\begin{lem} Let us denote by $\widetilde{\cl H}$ the elements of $\cl H$ which are of the form
\eqref{by} {\it with} $\phi(z)=\langle\lambda,z\rangle$, $\lambda\in E'$.
Then the elements of $\widetilde{\cl H}$ are exposed points with exposing
hyperplane $\lambda$ and such that $\Lambda(s\lambda)<+\infty$ for some $s>1$. Moreover for every open set $A\subset E$ we  have
\begin{equation}\label{hprime}
\inf_{A\cap \cl E}I=\inf_{A}I
\end{equation}
where $\cl E$ is as in the statement of Theorem \ref{egb}.
\end{lem}
\proof Remember, from \eqref{attained2}, that if $x\in\widetilde{\cl H}$ and $x=\int_E\langle\lambda,z\rangle z\,d\mu(z)$ then the supremum in \eqref{attained} is attained at
$\th=\lambda$ and that $I(x)=\Lambda(\lambda)=\langle\lambda,x\rangle-\Lambda(\lambda)$. Let us prove
that such an $x$ is an exposed point with exposing hyperplane $\lambda$.

Let $y\in\cl H$, $y\not=x$,  be such that $y=\int_Ez\phi(z)\,d\mu(z)$ for some $\phi\in \overline{E}'_\mu$.
Then
$$
\dlines{
I(y)=\frac 12\int_E\phi(z)^2\, d\mu(z)>
\frac 12\int_E\phi(z)^2\, d\mu(z)-\frac 12\underbrace{\int_E(\langle\lambda,z\rangle-\phi(z))^2\,d\mu(z)}_{>0}=\cr
=\int_E\langle\lambda,z\rangle\phi(z)\,d\mu(z)-
\frac 12\int_E\langle\lambda,z\rangle^2\,d\mu(z)=\langle \lambda,y\rangle-\Lambda(\lambda)\ .\cr
}
$$
Hence, as $I(x)=\langle \lambda,x\rangle-\Lambda(\lambda)$,
$$
I(y)-I(x)>\langle \lambda,y\rangle-\Lambda(\lambda)-(\langle \lambda,x\rangle-\Lambda(\lambda))=\langle \lambda,y-x\rangle
$$
so that every $x\in \widetilde{\cl H}$ is an exposed point. The condition $\Lambda(s\lambda)<+\infty$ for some $s>1$ is obviously satisfied as
$\Lambda(s\lambda)=s^2\Lambda(\lambda)$.

In order to prove \eqref{hprime} we need only to show that $\inf_{A\cap \cl E}I\le \inf_{A}I$. As $I\equiv +\infty$ outside $\cl H$,
$$
\inf_{A}I= \inf_{A\cap\cl H}I\ .
$$
Let us assume first that $\mathop{\rm supp} \mu=E$, so that $\cl H$ is dense in $E$. As $\widetilde{\cl H}$ is dense in $\cl H$, for every $x\in {A}\cap\cl H$ there exists a sequence $(x_n)_n\subset \widetilde{\cl H}$
converging to $x$ in the topology of $\cl H$, hence such that
$I(x_n)\to I(x)$. As the topology of $\cl H$ is stronger than the
topology of $E$, $(x_n)_n$ also converges to $x$ in
$E$, whence $x_n\in {A}$ for $n$ large (recall that $A$ is an open set). Therefore, for every $x\in\cl
H\cap A$ and $ n >0$ there exists $y\in\widetilde{\cl H}\cap A$ such that
$I(y)\le I(x)+ \frac 1n $. As $\widetilde{\cl H}\subset\cl E$, we can conclude as
$$
\inf_{A}I= \inf_{A\cap\cl H}I=\inf_{A\cap\widetilde{\cl H}}I\ge \inf_{A\cap \cl E}I\ .
$$
If $\mathop{\rm supp} \mu\not=E$ then $\widetilde E:=\mathop{\rm supp}\mu$ is a proper closed subspace of $E$ and it is immediate that $I\equiv+\infty$ on
$E\setminus \widetilde E$. Actually if $x\in E\setminus \widetilde E$, by the Hahn-Banach theorem there exists $\lambda\in E'$
such that $\langle\lambda,z\rangle=0$ for every $z\in \widetilde E$ and $\langle\lambda,x\rangle\not=0$. Hence
$$
I(x)=\sup_{\th\in E'}\bigl(\langle\th,x\rangle-\Lambda(\th)\bigr)\ge
\sup_{t\in \R}\bigl(t\langle\lambda,x\rangle-\Lambda(t\lambda)\bigr)=\sup_{t\in \R}t\langle\lambda,x\rangle=+\infty
$$
Moreover, considering $\mu$ as a Gaussian probability on $\widetilde E$, $\cl H$ is dense in $\widetilde E$
(in the topology of $\widetilde E$) and by a repetition of the previous argument we have
$$
\inf_{A}I= \inf_{A\cap\widetilde E}I=\inf_{A\cap\widetilde E\cap\cl H}I=\inf_{A\cap\widetilde E\cap\widetilde{\cl H}}I\ge
\inf_{A\cap\widetilde E\cap \cl E}I=\inf_{A\cap\cl E}I\ .
$$
\finedim
Remark that the argument above provides a proof of the LDP of Theorem \ref{mt2} which is new also in the classical case,
i.e. when $X_n\equiv X$.
At least to the author's knowledge.

\section{Second proof: exponential approximation}\label{exp-appr} The second proof of Theorem \ref{mt2}
is based on the notion of exponential equivalence (see \cite{DZ10}
\S4.2.2). Theorem 4.2.13 there states that if, for every $\delta>0$,
\begin{equation}\label{ee}
\limsup_{n\to 0} \frac1{g_n}\log\P(\vert g_n^{-1/2}\,Y_n
-g_n^{-1/2}\,X_n\vert\ge\delta)=-\infty
\end{equation}
then if $(g_n^{-1/2}\,Y_n)_n $ satisfies  a Large Deviation Principle
with speed $g$ and rate function $I$, the same is true for $(g_n^{-1/2}X_n)_n $. We shall apply this criterion with $Y_n\equiv X$.
\smallskip

By Skorokhod representation theorem, we can assume that the r.v.'s in the statement of Theorem \ref{mt2} are
defined on the same probability space and that $X_n\to X$ a.s. Hence
\begin{equation}\label{just}
\limn\frac 1{g_n}\log\E[\e^{g_n^{1/2}\langle \th ,X_n-X\rangle}]=\limn \frac12\Var(\langle \th,X_n-X\rangle)=0\ .
\end{equation}
By Theorem \ref{mt} the sequence $(g_n^{-1/2}(X_n-X))_n$ is
exponentially tight at speed $g$ and by \eqref{just} it satisfies the relation \eqref{eg.177} with
$\Lambda\equiv 0$. We can therefore apply the upper bound \eqref{eg.up0} of Theorem \ref{egb} with respect to the convex
conjugate $I$ of $\Lambda\equiv 0$, i.e. $I(0)=0$, $I(x)=+\infty$ for every $x\not=0$ and to the closed set $F=\{|x|\ge\delta\}$ in
order to obtain
\begin{equation}\label{ee2}
\limsup_{n\to \infty} \frac 1{g_n} \log\P(\vert g_n^{-1/2}\,X_n
-g_n^{-1/2}\,X\vert\ge\delta)=-\inf_FI=-\infty\ .
\end{equation}
Hence $(g_n^{-1/2}X_n)_n$ satisfies a Large Deviation principle with the same rate function as the sequence $(g_n^{-1/2}X)_n$ and
the result follows from the classical Large Deviation results of \cite{dvIII}.
\qed
It does not seem immediate to verify the exponential equivalence relation \eqref{ee} under the hypotheses of
Theorem \ref{mt2} without taking advantage of Theorems \ref{mt} and \ref{egb}.
\section{Examples and concluding remarks}\label{concl}%
\begin{example}\rm Let $B$ be a real Brownian motion and, for every integer $n$ and $t\le T$,
$$
X^{(n)}_t=\sqrt{2}\int_0^t\sin (ns) \, dB_s\ .
$$
We have $X^{(n)}_t=W^{(n)}_{A_n(t)}$ where $A_n(t)=2\int_0^t\sin^2(ns)\,ds$ and $W^{(n)}$ is another Brownian motion. As $A_n(t)\to t$ uniformly as
$ n \to\infty$, the sequence of processes $(X^{(n)})_n$ converges in law to the Wiener
measure. Hence $(X^{(n)})_n$ is a Gaussian tight sequence taking values in the Banach space $E=\cl C([0,T],\R)$ of continuous paths endowed
with the topology of uniform convergence and, by Theorem \ref{mt2}, $(g_n^{-1/2}X^{(n)})_n $ satisfies a Large Deviation
Principle at speed $g$ with the same rate function as in Schilder's theorem
(\cite{schilder}) for the Brownian motion.
\end{example}
\begin{example}\rm Let again $E=\cl C([0,T],\R)$ and let $X_n$ be a fractional Brownian motion with parameter $H_n$ with $\limn H_n=H$,
$0<H< 1$. It is well known that $X_n\to_{n\to\infty}X$, where $X$ is a fractional Brownian motion with Hurst index $H$, the convergence being
in law. This follows easily: the convergence of the covariance functions implies the convergence of the finite dimensional distributions and
the equality
$$
\E[(X_n(t)-X_n(s))^{2p}]=\frac {(2p)!}{2^pp!}\,|t-s|^{2pH_n}
$$
implies the tightness of the sequence thanks to Billingley's criterion (\cite{bill-cpm}, Theorem 12.3). Hence Theorem \ref{mt2}
implies immediately
that, for any speed function $g$, the sequence $(g_n^{-1/2}X_n)_n$ enjoys a Large Deviation estimate with
respect to the rate function given by the RKHS of the fractional Brownian motion with Hurst parameter equal to $H$, as can be found in
\cite{decreusefond} e.g.
\end{example}
It is well known that, if $E$ is finite dimensional, \eqref{eg.177} with the assumption that the limit $\Lambda$ is finite in a neighborhood of the origin implies
exponential tightness, even without Gaussianity. The author has been wondering whether  \eqref{eg.177} might imply exponential tightness
for Gaussian sequences in infinite dimensions.  Here we produce a counterexample.
\begin{example}\label{czero}\rm
Let us denote, as usual, by $c_0$ the space of real sequences $x=(x_n)_n$ such that $\limn x_n=0$ which is a separable Banach space with respect to the sup norm.

Let $(\xi_n)_n$ be a sequence of i.i.d. $N(0,1)$-distributed r.v.'s and
$$
X_n=(\xi_1,\xi_2,\dots,\xi_n,0,0,\dots)
$$
be a corresponding $c_0$-valued sequence of r.v.'s. We know that $(X_n)_n$ does not converge, hence it is not tight in $c_0$, because
this would require $\limn \xi_n=0$ a.s.

We prove first that $(\frac 1{\sqrt{n}}\, X_n)_n$ is not exponentially tight at speed $g_n=n$, i.e. that,
for $R>0$ fixed, there exist no compact
set $K=K_R\subset c_0$ such that, for $n\ge n_0$,
$$
\P\Bigl(\frac 1{\sqrt{n}}\, X_n\not\in K\Bigr)\le \e^{-nR}\ .
$$
Actually we prove that, for every compact set $K\subset c_0$,
\begin{equation}\label{no-et}
\liminf_{n\to\infty}\frac 1n\,\log \P\Bigl(\frac 1{\sqrt{n}}\, X_n\not\in K\Bigr)=0\ .
\end{equation}
It is easy to see that every such compact set $K$ is contained in a set of the form $\{x, |x_k|\le a_k \mbox{ for every }k=1,2,\dots\}$, where $(a_k)_k$ is a positive sequence converging to $0$ as $k\to +\infty$.
We have, for every $k_0>0$ and $n\ge k_0$,
$$
\dlines{
\P\Bigl(\frac 1{\sqrt{n}}\, X_n\not\in K\Bigr)\ge
\P\bigl(|\xi_k|> a_k\sqrt{n} \hbox{ for at least a } k=1,\dots,n\bigr)\ge\cr
\ge \P(|\xi_{k_0}|> a_{k_0}\sqrt{n}\,)
=2\bigl(1-\Phi(a_{k_0}\sqrt{n}\ )\bigr)\cr
}
$$
where by $\Phi$ we denote the $N(0,1)$ partition function. Now, thanks to classical estimates on the behavior of $\Phi$ at infinity,
$$
\liminf_{n\to\infty}\frac 1n\,\log \Bigl(2\bigl(1-\Phi(a_{k_0}\sqrt{n}\ )\bigr)\Bigr)\ge -\frac {a_{k_0}^2}2
$$
hence, $k_0$ being arbitrary and $(a_k)_k$ infinitesimal, we have \paref{no-et}.

However the limit \eqref{eg.177} exists in this case. Recall that here $E'=\ell_1$, the space of the sequences $\th=(\th_n)_n$ such that
$\sum_{n=1}^\infty|\th_n|<+\infty$.
As the r.v.
$$
\frac 1{\sqrt{n}}\, \langle \th,X_n\rangle=\frac 1{\sqrt{n}}\, \sum_{k=1}^n\th_k\xi_k
$$
is Gaussian, we have
$$
\log\E[\e^{n^{1/2}\langle \th,X_n\rangle}]=\frac n2\,\Var(\langle \th,X_n\rangle)=
\frac n2\,\sum_{k=1}^n\th_k^2
$$
hence the limit in \eqref{eg.177} exists and takes the value $\frac12\,\sum_{k=1}^\infty|\th_k|^2<+\infty$ for every $\th\in E'$.
\end{example}
\begin{example}\label{counter}\rm
Let $(\xi_n)_n$ be a sequence of i.i.d. $N(0,1)$-distributed r.v.'s and let
$$
X_n=(a_1\xi_1,\dots,a_n\xi_n,0,0\dots)
$$
be a $c_0$-valued r.v. where $(a_n)_n$ is a positive sequence converging to $0$. In this example we show that for some choice of the sequence $(a_n)_n$ the sequence $(X_n)_n$ is not tight but is exponentially tight, at least for some choices of the speed function.

Let us investigate the tightness of $(X_n)_n$.
Let $K=\{x\in c_0, |x_k|\le b_k\}$, where $(b_k)_k$ is a positive sequence converging to $0$, be a compact set of $c_0$. We have
$$
\P(X_n\in K^c)=1-\P(X_n\in K)=1-\prod_{k=1}^n
\P\Bigl(|\xi_k|\le \frac{b_k}{a_k}\Bigr)\ .
$$
The quantity $\P(|\xi_k|\le \frac{b_k}{a_k})$ appearing in the infinite product is always $\le 1$, hence it is well known that the infinite product converges to a number that is $>0$ if and only if the series
$$
\sum_{k=1}^\infty\Bigl(\P(|\xi_k|\le \frac{b_k}{a_k})-1\Bigr)
$$
is convergent. But  $\P(|\xi_k|\le \tfrac{b_k}{a_k})=2F(\tfrac{b_k}{a_k})-1$, with $F$ denoting the partition function of the standard Gaussian distribution. But recalling that, for large $x$,
$$
1-F(x)\ge \frac 1{2\pi}\,\Bigl(x+\frac 1x\Bigr)^{-1}\e^{-x^2/2}\ge \sqrt{\frac 2\pi}\,\frac1x\,\e^{-x^2/2}
$$
we have
$$
1-\P\Bigl(|\xi_k|\le \frac{b_k}{a_k}\Bigr)\ge
const\cdot\frac 1{b_k\sqrt{\log(k+1)}}\,\frac 1{(k+1)^{b^2_k/2}}\
$$
and this is the general term of a divergent series for every sequence $(b_k)_k$ such that  $b_k\to_{k\to 0}0$. Hence we have $\sup_{n\ge1}\P(X_n\in K^c)=1$ and $(X_n)_n$ is not tight.

Let us turn now to the exponential tightness of this sequence. Let $K=\{x\in c_0,\, |x_k|\le b_k\}$ be a compact set of $c_0$ and $R>0$ as above. We must find a sequence $(b_k)_k$ such that, for every $n$,
$$
\P(g_n^{-1/2}X_n\in K^c)\le \e^{-Rg_n}\ .
$$
Now
$$
\dlines{
\P(g_n^{-1/2}X_n\in K^c)=1-\P(g_n^{-1/2}X_n\in K)=1-\prod_{k=1}^n\P\Bigl(|\xi_k|\le g_n^{1/2}\,\frac{b_k}{a_k}\Bigr)=\cr
=1-\exp\Bigl\{\sum_{k=1}^n\log\P\Bigl(|\xi_k|\le g_n^{1/2}\,\frac{b_k}{a_k}\Bigr)\Bigr\}\ .
}
$$
As $1-\e^x\le -x$,
$$
\dlines{
\dots\le-\sum_{k=1}^n\log\P\Bigl(|\xi_k|\le g_n^{1/2}\,\frac{b_k}{a_k}\Bigr)=-\sum_{k=1}^n\log\Bigl(2F\Bigl(g_n^{1/2}\,\frac{b_k}{a_k}\Bigr)-1\Bigr)=\cr
=-\sum_{k=1}^n\log\Bigl(1-2\Bigl(1-F\Bigl(g_n^{1/2}\,\frac{b_k}{a_k}\Bigr)\Bigr)\Bigr)\ .\cr
}
$$
As $\log (1-x)\ge -x-\frac12\,x^ 2$,
\begin{equation}\label{eq-counter2}
\dots\le2\sum_{k=1}^n\Bigl(1-F\Bigl(g_n^{1/2}\,\frac{b_k}{a_k}
\Bigr)\Bigr)+\frac 12\,\Bigl(2\sum_{k=1}^n\Bigl(1-F\Bigl(g_n^{1/2}\,\frac{b_k}{a_k}
\Bigr)\Bigr)\Bigr)^2\ .
\end{equation}
With the well known inequality $1-F(x)\le \frac1{\sqrt{2\pi}}\,\e^{-x^2/2}$, $x\ge 0$,
$$
\sum_{k=1}^n\Bigl(1-F\Bigl(g_n^{1/2}\,\frac{b_k}{a_k}
\Bigr)\le\frac1{\sqrt{2\pi}}\,\sum_{k=1}^n\e^{-g_n\tfrac {b^2_k}{2a_k^2}}\ .
$$
{\it If $g_n\ge M\log n$, $n\ge 2$, for some $M>0$}, let us choose $b_k^2=2(R+\frac 1M)a_k^2$. We have
$$
\sum_{k=1}^n\Bigl(1-F\Bigl(g_n^{1/2}\,\frac{b_k}{a_k}
\Bigr)\le\frac1{\sqrt{2\pi}}\, \frac 1n\sum_{k=1}^n\e^{-R g_n}
=\frac1{\sqrt{2\pi}}\,\e^{-Rg_n}\ .
$$
Going back to \eqref{eq-counter2}, we find finally
$$
\P(g_n^{-1/2}X_n\in K^c)\le const\cdot\e^{-Rg_n}\ .
$$
so that {\it $(X_n)_n$ is exponentially tight, for every speed function increasing to $+\infty$ faster than $n\mapsto M\log n$}.

For this class of speed functions we can deduce that the sequence $(g_n^{-1/2}X_n)_n$ enjoys a Large Deviations Principle at
speed $g$. Actually the dual of $c_0$ is the space of summable sequences $\ell_1$ and, for every $\th\in\ell_1$ we have,
$$
\limn\frac 1{g_n}\,\log\E(\e^{g_n\langle \th,g_n^{-1/2}X_n\rangle})=\limn
\frac 12\,\Var(\langle \th,X_n\rangle)=\frac 12\,\sum_{k=1}^\infty\th_k^2a_k^2:=\Lambda(\th)
$$
and we find easily for the convex conjugate, $I$, of $\Lambda$, i.e.
$$
I(x)=\sup_{\th\in\ell_1}\bigl(\langle \th,x\rangle-\Lambda(\th)\bigr)\ ,
$$
that
$$
I(x)=\sum_{k=1}^\infty\frac {x_k^2}{a_k^2}
$$
with the understanding $I(x)=+\infty$ if the sum is not convergent.

It is now easy to complete this example by showing that actually the sequence $(g_n^{-1/2}X_n)_n$ enjoys a
Large Deviations Principle at speed $g$ and with respect to the rate function $I$. As the sequence $(a_k)_k$
tends to $0$ at infinity, it is immediate that the level sets of $I$ are compact in $c_0$. Moreover Theorem
\ref{egb} guarantees that the upper bound holds with respect to $I$ (this is equation \eqref{eg.up0}) and also
the partial lower bound \eqref{eg.low0}. In order to conclude we must only prove that the exposed points for $I$
are dense in $c_0$. Let us consider the set $\widetilde{\cl H}\subset c_0$ of the sequences $x=(x_k)_k$ that are
$=0$ for $k\ge N$ for some $N>0$. Of course $\widetilde{\cl H}$ is dense in $c_0$. The elements of $\widetilde{\cl H}$
are exposed points. Actually, let $\eta\in \ell_1$ be given by
$$
\eta_k=\frac {x_k}{a_k^2}\ \raise2pt\hbox{,} \qquad k\le N
$$
and $\eta_k=0$ for $k>N$. Then if $x\in\widetilde{\cl H}$ we have for every $y\in c_0$,
$y\not=x$,
$$
\dlines{
I(y)-I(x)-\langle\eta,y-x\rangle=\frac 12\sum_{k=1}^\infty \frac {y_k^2}{a_k^2}-\frac 12\sum_{k=1}^\infty \frac {x_k^2}{a_k^2}-\sum_{k=1}^\infty \frac {x_k}{a_k^2}\,(y_k-x_k)=\cr
=\frac 12\sum_{k=1}^\infty \frac {y_k^2}{a_k^2}+\frac 12\sum_{k=1}^\infty \frac {x_k^2}{a_k^2}-\sum_{k=1}^\infty \frac {y_kx_k}{a_k^2}=\frac 12\sum_{k=1}^\infty\frac 1{a_k^2}\,(y_k-x_k)^2>0
}$$
hence $\eta$ is an exposing hyperplane and we are allowed to conclude.
\end{example}
\bigskip

\bibliography{bibbase}
\bibliographystyle{amsplain}
\end{document}